\documentclass[twoside, 11pt]{article}
\usepackage[section]{algorithm}
\usepackage{amsthm,amsmath,amssymb,amsmath,amsfonts,graphicx,fancyhdr,algorithmic,longtable,listings,color}
\usepackage{bibcheck}
\usepackage{graphicx}

\makeatletter
\newcounter{algorithmbis}
\renewcommand{\thealgorithmbis}{\thesection.\arabic{algorithmbis}}
\def\algorithmbis{\@ifnextchar[{\@algorithmbisa}{\@algorithmbisb}}
\def\@algorithmbisa[#1]{%
  \refstepcounter{algorithmbis}
  \trivlist
  \leftmargin\z@
  \itemindent\z@
  \labelsep\z@
  \item[\parbox{\textwidth}{%
    \hrule
    \hrule
    \noindent\strut\textbf{Algorithm \thealgorithmbis} #1
    \hrule
  }]\hfil\vskip0em%
}
\def\@algorithmbisb{\@algorithmbisa[]}

\makeatother

\def\enddemo{\qed \endtrivlist}
\expandafter\let\csname enddemo*\endcsname=\enddemo

\def\qedsymbol{\ifmmode\bgroup\else$\bgroup\aftergroup$\fi
  \vcenter{\hrule\hbox{\vrule height.6em\kern.6em\vrule}\hrule}\egroup}
\def\qed{\ifmmode\else\unskip\nobreak\fi\quad\qedsymbol}

\newtheorem{exm}{Example}[section]

\newcounter{tabela}[section]

\def\B{{\cal{B}}}

\def\R{{\cal{R}}}

  \font\sst=cmtt8 
  \font\ssi=cmti8 
\font\srs=cmr6

\def\ddj{d\kern -0.70em\lower 0.22em\hbox{\B\ }\kern -0.02em}

\def\dj{d\kern-0.4em\char"16\kern-0.1em}
\def\Dj{\hbox{\raise0.3ex\hbox{-}\kern-0.4em D}}

\title{\bf Single-facility Weber Location Problem \\ based on the Lift Metric}

\author{\frenchspacing
\bf Predrag S. Stanimirovi\' c$^{1}$\footnote{Corresponding author}, \bf  Marija \'Ciri\'c$^{2}$ \\
${}^{\srs 1,2}${\ssi University of Ni\v{s}, Faculty of Sciences and Mathematics,}\\
{\ssi Vi\v{s}egradska 33, 18000 Ni\v{s}, Serbia.}
\footnote{The authors gratefully acknowledge support from the research project 144011 of the Serbian Ministry of Science.} \\
{\ssi {\ssi E-mail}:} $^{\srs 1}${\sst  pecko@pmf.ni.ac.rs},\ $^{\srs 2}${\sst marijamath@yahoo.com}}
\date{}

\begin{document}
\maketitle

\begin{abstract}
The continuous single-facility min-sum Weber location problem based upon the lift metric is investigated.
An effective algorithm is developed for its solution.
Implementation for both the discrete and continuous location problems is developed in the programming package {\it Mathematica}.
\begin{description}
\frenchspacing \itemsep=-1pt
\item[] AMS Subj.  Class.: 90B80, 90B85 \/
\item[] Key words: Continuous location problem, Weber problem, Lift metric, {\it Mathematica}.
\end{description}
\end{abstract}

\section {Introduction\/}\setcounter{equation}{0}

Location problems represent a special class of optimization tasks, where coordinate of locations and distances between them are main parameters.
In the general case, the task of location problem is to
define positions of some new facilities from the actual space in which are already placed some other relevant objects (points).
New facilities are centers that provide services and called {\it suppliers}.
Existing facilities are the service users or clients, and called {\it customers}.
Location problems occur frequently in real life.
Many systems in the public and private sectors are characterized by facilities that provide
homogeneous services at their locations to a given set of fixed points or customers. Examples
of such facilities include warehouse location, positioning a computer and communication units, locating hospitals, police stations,
locating fire stations in a city, locating base stations in wireless networks.

\smallskip
Different classifications of the location problems are known.
The classification scheme from \cite{Hamacher1} assumes five positions.

In the present article we pay attention to the selection of the distance function as the characterization criterion of the location problem.
The distance between two points is the length of the shortest path connecting them.
The metric by which the (generalized) distance between two points is measured may be different in various instances \cite{Chen}.
In the calculating of distance between two points,
the most common distance metrics in a continuous space are those known as the class of $l_p$ distance metrics,
primarily rectangular ($l_1$), Euclidean ($l_2$) and Chebyshev ($l_\infty$) metric.
Detailed explanation of various metrics one can find in Dictionary of distances \cite{[1]}.
Many factors affect on the process of metrics choosing.
The most important factor is the nature of the problem. For example, if it is possible to move rectilinearly between
two points, the distance between them is exactly given by the Euclidean (or straight-line distance) metric.
On the other hand, in the cities where streets intersect
under the right angle mainly, the distance between two points will be the
best approximated using the rectangular metric (also known as the Manhattan, "city block" distance, the right-angle distance metric
or taxicab distance).
Measures of distances in chess are a characteristic example.
The distance between squares on the chessboard for rooks is measured in Manhattan distance;
kings and queens use Chebyshev distance, and bishops use the Manhattan distance.

\smallskip
We emphasize the next main contribution of our paper.

\smallskip
\noindent The Weber location problem (also called the Fermat-Weber problem) is a basic model in the location
theory which has received significant attention in the scientific literature.
For a detailed review see, for example, \cite{Wesolowsky}.
The paper \cite{[7]} investigated a reformulation of the unconstrained form of the classical Weber
problem into an unconstrained minimum norm problem.
The classical Weber problem is established with the Euclidean norm underlying in the definition of the distance function.
But, other measures, principally $l_p$ norms, also play an important role in the theory and
practice of location problems. The norms are arbitrary, in general.
The most popular method to solve the Weber problem problem with Euclidean distances
is given by a one-point iterative procedure which was first proposed by Weiszfeld \cite{Weiszfeld}.
The procedure is readily generalized to $l_p$ distances (see, for example, \cite{Love}, Ch. 2).
Solution of the continuous Weber problem in $l_1$ distance is described in \cite{[2]}.
The three-dimensional Fermat-Weber facility location problem with Tchebychev
distance is investigated in \cite{[8]}.
The Weber location problem with squared Euclidean distances is considered in \cite{[2]};
the same problem under the assumption that the weights are selected from a given set of intervals at any point,
is studied in \cite{Drezner2}.

\smallskip
The $l_p$ norms have received the most attention from location analysts.
But, many other types of distances have been exploited in the facility location problem.
A review of exploited metrics is presented in \cite{[2]}:

- central metrics \cite{Perreur}, 

- distance functions based on altered norms \cite{Love1,Love2},

- weighted one-infinity norms \cite{Ward}, 

- mixed norms \cite{Hansen}, 

- block and round norms \cite{Thisse}, 

- mixed gauges \cite{Durier}, 

- asymptotic distances \cite{Hodgson}, 

- weighted sums of order $p$ \cite{Brimberg1,Uster}. 

\smallskip
In the present article we solve the Weber problem in the plane, under the assumption that the distance is measured by the lift metric.



\smallskip
The paper is organized as follows.
Some basic definitions and algorithms are restated in the second section.
In the third section we present an effective algorithm for the solution of the single-facility continuous Weber problem,
assuming that distances are measured by the lift metric.

\section {Preliminaries\/}\setcounter{equation}{0}

The {\it lift metric} or the {\it raspberry picker metric} in the plane $\R^2$ is defined by
\begin{equation}\label{Lift}
L(A,B)=\left\{\begin{array}{ll}
|x_1^A-x_1^B|,  &x_2^A=x_2^B\\ \smallskip
|x_1^A|+|x_2^A-x_2^B|+|x_1^B|, & x_2^A\neq x_2^B
\end{array}
\right.
\end{equation}
where $A(x_1^A,x_2^A)$ and $B(x_1^B,x_2^B)$ are given points.
It can be defined as the minimum Euclidean length of all admissible connecting curves
between two given points, where a curve is called admissible if it consists of only segments of straight lines parallel to $x$-axis, and
of segments of $y$-axis (see, for example \cite{[1]}).
Therefore, under the assumption $x_2^A\neq x_2^B$ the
distance between two points $A$ and $B$ in the lift metric equals the sum of lengths $AA'$, $A'B'$ and $B'B$, where $A'$
and $B'$ are orthogonal projections of the points $A$ and $B$ to the $y$-axis, respectively (Figure 1,Left).
In the opposite case, $x_2^A=x_2^B$, the distance between $A$ and $B$ is simply the length of the segment $AB$ (Figure 1,Right).

\begin{figure}[hbtp]
\centering $\includegraphics* [scale = 0.6] {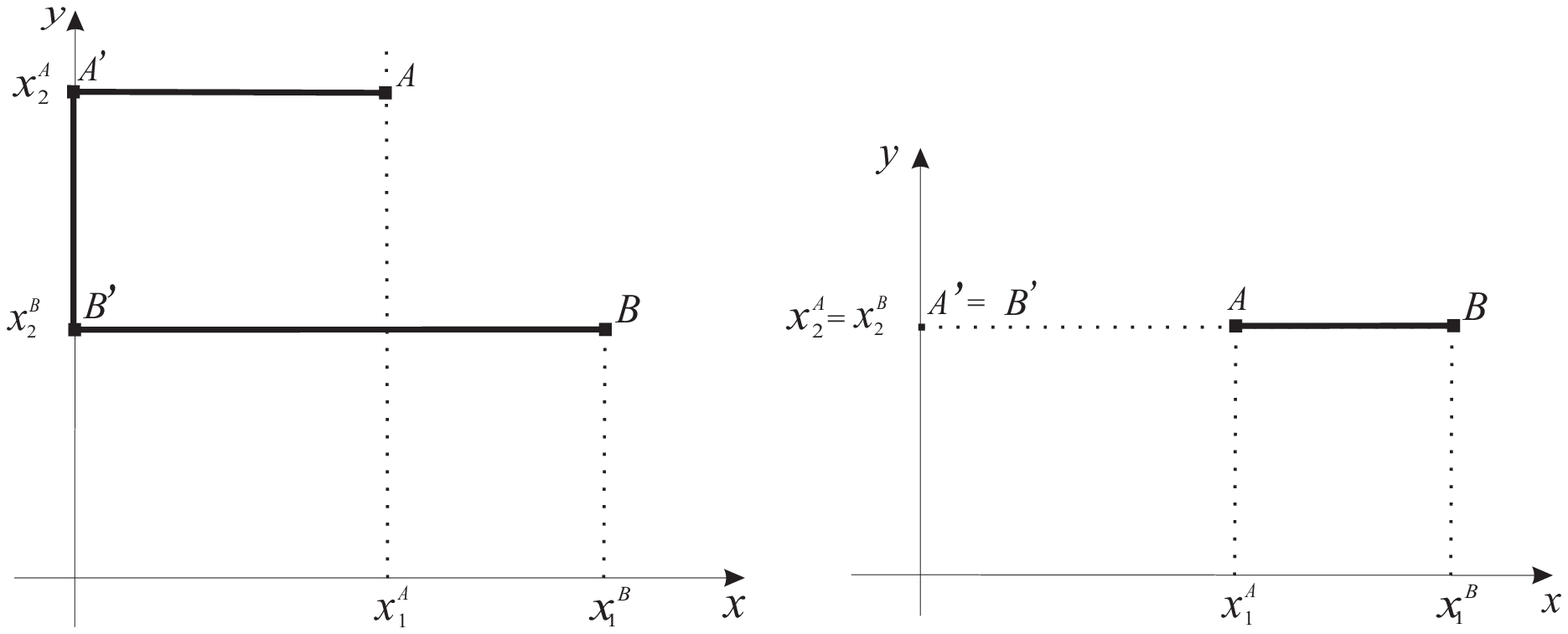}
$ \center{\small Figure 1. Left) The case $x_2^A\neq x_2^B$ \hskip 0.5cm  Right) The case $x_2^A= x_2^B$}
\end{figure}

This distance is appropriate for usage in cities which have one main street (corresponding to the $y$-axis),
and the other side streets are normal to it (Figure 2). We are observed that
in the main city of Zakynthos island in Greek-Zakynthos, the streets are deployed on this way.
Similar situation also occurs in tier buildings where the lift (in the role of $y$-axis) connects tiers.

The 2-dimensional continuous Weber location problem can be briefly restated as follows (see, for example \cite{Drezner,Wesolowsky}).
Let $m$ demand centers $A_1,\ldots ,A_m$ be given in the plane $\R^2$ (locations of given customers),
where $A_i(a_1^i,a_2^i),\mbox{  }i=1,\ldots ,m.$
It is necessarily to find a new point $X(x_1,x_2)\in \R^2$ which has minimal sum of weighted distances with respect to given points.
Therefore, one needs to solve the unconstrained optimization problem (single-facility min-sum problem),
where it is necessary to minimize the sum
\begin{equation}\label{(1)}
\min\limits _X \ f(X)=\sum_{i=1}^m w_i\cdot d(A_i,X).
\end{equation}
The real quantity $w_i$ is a positive weighted coefficient of the point $A_i$.
Essentially, the weight $w_i$ converts the distance $d(A_i,B_k)$ 
into a cost of serving the demand of customer $i$ considerate to $k$th offered facility location.

\smallskip
For the sake of completeness, we restate well-known method which gives solution of the Weber problem \eqref{(1)}
in the case when the underlying distance function is defined by the $l_1$ metric, which means
$$l_1(A_i,X)=|x_1-a_1^i|+|x_2-a_2^i|.$$
Then the goal function $f(X)$ divides in two sums:
$$f_1(x_1)=\sum_{i=1}^m w_i |x_1-a_1^i|,\quad f_2(x_2)=\sum_{i=1}^m w_i |x_2-a_2^i|. $$
Therefore, the initial optimization problem splits into two independent optimization tasks.
It is sufficient to consider the optimization problem
\begin{equation}\label{nova}
\min\limits _x \ f(x)=\sum_{i=1}^m w_i |x-a^i|.
\end{equation}
We restate well-known algorithm for solving the problem \eqref{nova} (one can find it in \cite{Mladen,[10]}).
The first step in this procedure is optional, but it can be used to accelerate the remaining two steps.

\begin{algorithmbis}[Solve the single-facility continuous Weber location problem.]\label{algWeb}
\begin{algorithmic}[1]
\REQUIRE Real quantities $a^1,\ldots ,a^m$.
  \STATE Step $I.$ For each subset of identical elements $a^{i_1}=a^{i_2}=\cdots =a^{i_j}$ perform the following activities:
put $w_{i_1}=w_{i_1}+w_{i_2}\cdots +w_{i_j}$, eliminate multiple
elements $a^{i_2},\ldots ,a^{i_j}$ as well as corresponding weights $w_{i_2},\ldots ,w_{i_j}$
and later perform appropriate shifting of the indices of residual elements $a^j$ and their weights $w_j$.
  \STATE Step $II.$
If Step $I$ is applied, denote by $q$ the cardinal number of different elements in the set $\{a^{1},\ldots ,a^{m}\}$;
otherwise, use $q=m$.
Sort coordinates $a^i$, $i=1,\ldots,q$ in non-descending order.
Let the sorted sequence of coordinates is
$$a^{1'}\leq a^{2'}\leq \cdots\leq a^{q'}.$$
Rearrange the weighting coefficients $\{w_1,\ldots ,w_q\}\to \{w_{1'},\ldots ,w_{q'}\}$ applying identical replacements on the weights.
For the sake of simplicity, let us denote partial sums of the array $\{w_{1'},\ldots ,w_{q'}\}$ by
\begin{equation}\label{S2}
S_2[0]=0,\ \ \ S_2[k]=\sum_{i=1}^k w_{i'},\ 1\leq k\leq q.
\end{equation}
  \STATE Step $III.$ There are two possibilities, denoted by $P_1$ and $P_2$.

\smallskip
$P_1$.  If the condition
\begin{equation}\label{prva1}
S_2[k^*-1]<\frac 12 S_2[q]<S_2[k^*]    
\end{equation}
is satisfied 
for some $k^*\in \{1,\ldots ,q\}$, then we end the algorithm without any possible solution.
Indeed, formal solution $x=a^{k^*}$ is eliminated according to assumption \eqref{suprotno} in this case.

\smallskip
$P_2$. If the condition
\begin{equation}\label{druga1}
\frac 12 S_2[q]=S_2[k^*]
\end{equation}
holds 
for some $k^*\in \{1,\ldots ,q\}$, then the solution is multiple, i.e. the
searched coordinate $x_2$ can to have any value from the interval $(a^{k^*},a^{k^*+1})$.
As agreed, we use the midpoint value $x=(a^{k^*}+a^{k^*+1})/2$.
\end{algorithmic}
\end{algorithmbis}

\section {Continuous Weber problem and lift metric\/}\setcounter{equation}{0}

In the sequel we solve the single-facility min-sum Weber problem \eqref{(1)} applying the lift metric \eqref{Lift}.
Therefore, the distance function is defined by
\begin{equation}\label{Lift1}
d(A_i,X)=L(A_i,X)\left\{\begin{array}{ll}
|x_1-a_1^i|,  &x_2=a_2^i\\ \smallskip
|x_1|+|x_2-a_2^i|+|a_1^i|, & x_2\neq a_2^i
\end{array}
\right. ,\ \ i=1,\ldots ,m.
\end{equation}

Two major steps (denoted as {\bf Step 1} and {\bf Step 2}) are separated in our algorithm, as in the following.

\smallskip
\noindent {\bf Step 1.} Generate the list $\mathfrak{X}$ of permissible solutions of the problem.
Its initial value is the empty set $\mathfrak{X}=\emptyset$.
Two different procedures are separated during the construction of the set $\mathfrak{X}$, in accordance with the definition \eqref{Lift1}.

\smallskip
{\bf Procedure 1.}
Let us consider the quotient set $S/\!\!\!=$ of the set $S=\{a_2^1,\ldots ,a_2^m\}$ into the equivalence classes
$S=\{S_1=[a_2^{i_1}],\ldots ,S_d=[a_2^{i_d}]\}$, where the class $S_j$ contains elements from $S$ whose values are $a_2^{i_j}$.
For each $j=1,\ldots ,d$ seek the second coordinate of the optimal point $X(x_1,x_2)$ in the form
\begin{equation}\label{Case1}
x_2\in S_j\Leftrightarrow x_2=a_2^{i_j}.
\end{equation}
So, as we know value of the coordinate $x_2$ ($x_2=a_2^{i_j}$), it is necessary to determine value of the coordinate $x_1$.
%
Denote  by $Q_j$ 
the set of indices corresponding to points whose second coordinates are contained in the set $S_j$.
According to \eqref{Lift1} and \eqref{Case1}, the objective function $f(X)$ consists of two separated sums
\begin{equation}\label{(2)}
f(X)= \sum_{i\notin Q_j}^m \!\!
w_i\left( |x_1|+ |x_2-a_2^i|+ |a_1^i|\right) +   \sum_{i\in Q_j}^m w_{i}|x_1\!-\!a_1^{i}|
\end{equation}
Taking into account $x_2=a_2^{i_j}$ and grouping the first term in the first sum with the second sum, we obtain
\begin{equation}\label{(3)}
f(X)=\sum_{i\notin Q_j}^m w_i\left( |a_2^{i_j}-a_2^i|+|a_1^i|\right)+ \sum_{i=1}^{m} w_{i}|x_1-a_1^{\beta(i)}|,
\end{equation}
where
\begin{equation}\label{novo}
a_1^{\beta(i)}=\left\{\begin{array}{rl}
 a_1^{i_p}, & \beta(i)=i_p\in I\\
 0, & \beta(i)\not\in I.
 \end{array}
 \right. 
\end{equation}
As the first sum in the expression \eqref{(3)} is constant, the
problem is reduced on determining the minimum of the function
\begin{equation}\label{(4)}
f_1(x_1)=\sum_{i=1}^{m} w_{i}|x_1-a_1^{\beta(i)}|.
\end{equation}

\noindent
We apply Algorithm \ref{algWeb} in adapted form for our specific situation \eqref{novo}, \eqref{(4)}.
That process consists of three major steps.

\smallskip
Step $I.$
For each subset of identical elements $a_1^{\beta(i_1)}=a_1^{\beta(i_2)}=\cdots =a_1^{\beta(i_j)}$ perform the following:
put $w_{i_1}=w_{i_1}+w_{i_2}+\cdots +w_{i_j}$, eliminate multiple
elements $a_1^{\beta(i_2)},\ldots ,a_1^{\beta(i_j)}$ as well as their weights $w_{i_2},\ldots ,w_{i_j}$
and then perform appropriate renumeration of the indices of the remainder elements $a_1^{\beta(j)}$ and their weights $w_j$.

\smallskip
Step $II.$
If Step $I$ is applied, denote by $p$ the cardinal number of different elements in the set $\{a_1^{\beta(1)}, \ldots ,a_1^{\beta(m)}\}$;
otherwise, use $p=m$.
Sort coordinates $a_1^{\beta(1)},\ldots ,a_1^{\beta(p)}$ in non-descending array.
Furthermore we suppose that the ordered sequence is
$$a_1^{1'}\leq a_1^{2'}\leq \cdots \leq a_1^{p'}.$$
Rearrange the corresponding weighting coefficients $w_{1},\ldots ,w_{p}$ analogously in the sequence
$$w_{1'},w_{2'}, \ldots ,w_{p'}.$$
For the sake of simplicity, let us denote partial sums of the array $\{w_{1'},\ldots ,w_{p'}\}$ by
$$S_1[0]=0,\ \ \ S_1[k]=\sum_{i=1}^k w_{i'},\ 1\leq k\leq p.$$

\smallskip
Step $III.$ There are two possible cases capable to produce permissible minimizers for $f_1$, denoted by $C_1$ and $C_2$.

\smallskip
$C_1.$ If the inequalities
\begin{equation}\label{prva}
S_1[k'-1]<\frac 12 S_1[p]<S_1[k'],\     
\end{equation}
are satisfied for some  $k'\in\{1,\ldots ,p\}$, then the searched coordinate is $x_1=a_1^{k'}.$
Later, we use $X(x_1,x_2)$ as the possible optimal point: $\mathfrak{X}=\mathfrak{X}\cup \{(a_1^{k'},a_2^{i_j})\}$.

\smallskip
$C_2.$  If the condition
\begin{equation}\label{druga}
\frac 12 S_1[p]=S_1[k']
\end{equation}
is satisfied for some $k'\in\{1,\ldots ,p\}$, then the solution is multiple, i.e.
the searched coordinate $x_1$ can to have any value from the interval $[a_1^{k'},a_1^{k'+1}]$.
In our implementation we use the value $x_1=(a_1^{k'}+a_1^{k'+1})/2.$
Thus, we found the additional possible solution of the starting problem \eqref{(1)},
which implies:
$\mathfrak{X}=\mathfrak{X}\cup \{((a_1^{k'}+a_1^{k'+1})/2,a_2^{i_1})\}$.

\smallskip
{\bf Procedure 2.}
Compute $x_2$ under the assumption
\begin{equation}\label{suprotno}
x_2\notin S \Leftrightarrow x_2\neq a_2^i\ \textrm{ for each } i\in \{1,\ldots,m\}
\end{equation}
(under the assumptions opposite with respect to \eqref{Case1}), the function $f(X)$ is reduced to
\begin{equation}\label{(5)}
f(X)=\sum_{i=1}^m w_i \left( |x_1|+ |x_2-a_2^i|+ |a_1^i|\right).
\end{equation}

\noindent It is necessary to minimize that function.
Since the third term $w_i|a_1^i|$ in the function $f(X)$ defined in \eqref{(5)} is constant,
one needs to minimize the next two objectives:
\begin{eqnarray}
\min\limits _{x_1}  \ f_1(x_1)\!\!\!\!&=&\!\!\!\! \sum_{i=1}^{m}w_i|x_1|\label{(6)}\\
\min\limits _{x_2} \ f_2(x_2)\!\!\!\!&=&\!\!\!\! \sum_{i=1}^m w_i|x_2-a_2^i|. \label{(7)}
\end{eqnarray}

Thus, solving the  problem \eqref{(5)} with two variables was reduced to
solving two independent tasks of unconstrained optimization \eqref{(6)} and \eqref{(7)} with one variable ($x_1$ and $x_2$, respectively).

\smallskip
Solution of the optimization problem \eqref{(6)} is evidently $x_1=0$.
The optimization problem \eqref{(7)} is the classical continuous Weber location model with underlying $l_1$ metric,
where only the additional assumption \eqref{suprotno} is imposed.
Therefore, in order to find optimal value for $x_2$ it suffices to apply Algorithm \ref{algWeb}
assuming that the input sequence is $a^1=a_1^1,\ldots ,a^m=a_2^m$ and taking into account conditions \eqref{suprotno}.

\smallskip
Step $I.$
For each subset of identical elements $a_2^{i_1}=a_2^{i_2}=\cdots =a_2^{i_j}$ perform the following activities:
put $w_{i_1}=w_{i_1}+w_{i_2}+\cdots +w_{i_j}$, eliminate multiple
elements $a_2^{i_2},\ldots ,a_2^{i_j}$ as well as corresponding weights $w_{i_2},\ldots ,w_{i_j}$
and later perform appropriate shifting of the indices of residual elements $a_2^j$ and their weights $w_j$.

\smallskip
Step $II.$
If Step $I$ is applied, denote by $q$ the cardinal number of different elements in the set $\{a_2^{1},\ldots ,a_2^{m}\}$;
otherwise, use $q=m$.
Sort coordinates $a_2^i$, $i=1,\ldots,q$ in non-descending order.
Let the sorted sequence of coordinates is
$$a_2^{1'}\leq a_2^{2'}\leq \cdots\leq a_2^{q'}.$$
Rearrange the weighting coefficients $\{w_1,\ldots ,w_q\}\to \{w_{1'},\ldots ,w_{q'}\}$ applying identical replacements on the weights.
Subsequently, generate the partial sums $S_2[i],i=0,\ldots ,q$ of the array $\{w_{1'},\ldots ,w_{q'}\}$ as in \eqref{S2}.

\smallskip
Step $III.$ There are two possibilities, denoted by $P_1$ and $P_2$.

\smallskip
$P_1$. If the condition \eqref{prva1}
is satisfied 
for some $k^*\in \{1,\ldots ,q\}$, then the algorithm is finished without any solution.
Indeed, the formal solution $x_2=a_2^{k^*}$ is eliminated according to assumption \eqref{suprotno}, actual for this case.

\smallskip
$P_2$. If the condition \eqref{druga1}
holds 
for some $k^*\in \{1,\ldots ,q\}$, then the solution is multiple, i.e. the
searched coordinate $x_2$ can to have any value from the interval $(a_2^{k^*},a_2^{k^*+1})$.
We use the midpoint value $x_2=(a_2^{k^*}+a_2^{k^*+1})/2$, so that the possible optimal point is $X(0,x_2)$.
Place the point $X$ at the end of the list $\mathfrak{X}$ by $\mathfrak{X}=\mathfrak{X}\cup \{(0,(a_2^{k^*}+a_2^{k^*+1})/2)\}$.

\smallskip
\smallskip
\noindent {\bf Step 2.}
Thus, we got one or more permissible solutions of the starting
problem \eqref{(1)}. For all obtained values $X$ from $\mathfrak{X}$ we
determine the values of the function $f(X)$ defined in \eqref{(5)}. Solution of the Weber
problem will be the point $X^*(x_1^*,x_2^*)$ for which the function $f(X)$ has a minimal value.
Actually in this step we are solving generated discrete location problem, where the set $\mathfrak{X}$ contains
in advance defined feasible locations of the supplier.

\smallskip
Let $X_1,\ldots ,X_r$ be $r$ locations on which it is possible to set a new desired object (supplier).
{\it The sum of weighted distances} from the permissible location $X_k$, $k\in \{1,\ldots ,r\}$ of the supplier to the customers is equal to
\begin{equation}\label{Ldisct}
W_k=\sum_{i=1}^m w_i\cdot L(A_i,X_k).
\end{equation}
The task is to determine the location $B_{k^*}$ for which the sum of weighted distances is minimal, i.e.
$$W_{k^*}=\min\ \{W_k\vert \ 1\leq k\leq r\}.$$

\smallskip
In accordance with the previous considerations, we state the following general algorithm.

\begin{algorithmbis}[Solution of the single-facility min-sum Weber problem in the lift metric.]\label{alg1}
\begin{algorithmic}[1]
\REQUIRE List $lp=\{(a_1^1,a_1^2),\ldots , (a_m^1,a_m^2)\}$ and the list of weights $lt=\{w_1,\ldots ,w_m\}$.
  \STATE {\bf Step 1:} Form the quotient set of $S=\{a_2^1,\ldots ,a_2^m\}$ in the form $\{S_1,\ldots ,S_d\}$, where
     each equivalence class $S_j$ contains identical elements from $S$ with the value $a_2^{i_j}$.
  \STATE {\bf Step 1:} Generate the list $\mathfrak{X}$ applying the procedure included into the possibilities $C_1,C_2$
  (included in {\bf Procedure 1.}) to all distinctive
  values $a_2^{i_j}$ of the set $S$, i.e. using $x_2= a_2^{i_j}$, $j=1,\ldots ,d$.
  \STATE {\bf Step 1:} Extend the list $\mathfrak{X}$ applying the method defined in the case $P_2$ (included in {\bf Procedure 2.}).
  \STATE {\bf Step 2:} Solve the discrete location problem using given locations $lp$, discrete set $\mathfrak{X}$ of possible solutions
  and the weights $lt$.

\end{algorithmic}
\end{algorithmbis}

\begin{footnotesize}
\begin{exm}\rm 
Solve Weber problem using the specified algorithm with the next data:
$$A_1(4,4),\mbox{  }w_1=4,\quad A_2(3,1),\mbox{  }w_2=1,\quad A_3(6,4),\mbox{  }w_3=2, \quad A_4(6,2),\mbox{  }w_4=3.$$
We have $S=\{4,1,4,2\}$. The quotient set of $S$ is defined as $S_1=[4],S_2=[1],S_3=[2]$.
Therefore, it is necessary to consider three possibilities for the cases $C_1$ and $C_2$.

\smallskip
{\bf 1.} Let be $x_2=4$. Then the function $f(X)$ has the
following  form:
$$
\aligned
f(x)=\ & w_2( |x_1|+|x_2-a_2^2|+|a_1^2|)+ w_4( |x_1|+|x_2-a_2^4|+|a_1^4|)+w_1|x_1-a_1^1|+w_3|x_1-a_1^3|.
\endaligned
$$
According to the constant value $x_2=4$ of the coordinate $x_2$, the function $f$
just depend on $x_1$, so we can consider the next function
$$
\aligned
f_1(x_1)&=w_1|x_1-a_1^1|+w_2|x_1-0|+w_3|x_1-a_1^3|+w_4|x_1-0|\\
&=\sum_{i=1}^4 w_{i}|x_1-a_1^{\beta(i)}|,
\endaligned
$$
where $a_1^{\beta(1)}=a_1^1=4$,\ $a_1^{\beta(3)}=a_1^3=6$,\ $a_1^{\beta(2)}=a_1^{\beta(4)}=0$.

\smallskip
Let us sort the coordinates $a_1^{\beta(i)}\to a_1^{i'}$ and rearrange corresponding weights $w_i\to w_{i'}$, using the same replacements:

\begin{center} \begin{tabular}{|l|l|l|l|l|}
\multicolumn{5}{r} {Table 1. }\\ \hline
 coordinates $(a_1^{i'})$ & 0 & 0 & 4 & 6 \\ \hline
 weights $(w_{i'})$ & 1 & 3 & 4 & 2\\ \hline
$k$ & 1 & 2 & 3 &4 \\ \hline
$S_1[k]=\sum_{i=1}^{k} w_{i'}$ & 1 & 4 & 8 & 10
\\\hline
\end{tabular}
\end{center}
\smallskip

\noindent
We firstly assume that Step $I$ is omitted.
According to $\frac 12\sum_{i=1}^4 w_{i'}=5$, the condition
$$S_1[k'-1]<\frac 12 S_1[4]<S_1[k']$$
is satisfied for $k'=3$, so $x_1=a_1^{3'}=4$.
Therefore, one possible solution is $X_1(4,4)$.

\smallskip
In the case when Step $I$ is applied, data from Table 1 reduce to

\begin{center} \begin{tabular}{|l|l|l|l|}
\multicolumn{4}{r} {Table 2. }\\ \hline
 coordinates $(a_1^{i'})$ & 0 & 4 & 6 \\ \hline
 weights $(w_{i'})$ & 4 & 4 & 2\\ \hline
$k$ & 1 & 2 & 3  \\ \hline
$S_1[k]=\sum_{i=1}^{k} w_{i'}$ & 4 & 8 & 10
\\\hline
\end{tabular}
\end{center}

Then conditions \eqref{prva} are satisfied for $k=2$, so that the same possible solution is generated.

\smallskip
The list of permissible solutions is now equal to $\mathfrak{X}=\{X_1\}$.

\medskip
{\bf 2.} In this case it is assumed $x_2=1$. Now the function $f_1(x_1)$ looks like:
$$f_1(x_1)=w_1|x_1-0|+w_2|x_1-3|+w_3|x_1-0|+w_4|x_1-0|.$$
On the similar procedure as in the case {\bf 1.} we get the table:

\begin{center}
\begin{tabular}{|l|l|l|l|l|}
\multicolumn{5}{r} {Table 3. }\\ \hline
coordinates $(a_1^{i'})$ & 0 & 0 & 0 & 3 \\ \hline
weights $(w_{i'})$ & 4 & 2 & 3 & 1\\ \hline
$k$ & 1 & 2 & 3 &4 \\ \hline
$S_1[k]=\sum_{i=1}^{k} w_{i'}$ & 4 & 6 & 9 & 10 \\\hline
\end{tabular}
\end{center}

\smallskip
\noindent Inequalities \eqref{prva}
are valid for $k'=2$, so $x_1=a_1^{2'}=0$, i.e.
we got the second possible solution $X_2(0,1)$.

\smallskip
In the case when Step $I$ is applied, Table 3 transforms to

\begin{center}
\begin{tabular}{|l|l|l|}
\multicolumn{3}{r} {Table 4. }\\ \hline
coordinates $(a_1^{i'})$ & 0 & 3 \\ \hline
weights $(w_{i'})$ & 9 & 1\\ \hline
$k$ & 1 & 2 \\ \hline
$S_1[k]=\sum_{i=1}^{k} w_{i'}$ & 9 & 10 \\ \hline
\end{tabular}
\end{center}

Condition \eqref{prva} holds for $k=1$, so that $X_2$ is again the second eventual solution.

\smallskip
We have $\mathfrak{X}=\{X_1,X_2\}$.

\medskip
{\bf 3.} Let us now start from the assumption $x_2=2$.
$$f_1(x_1)=w_1|x_1-0|+w_2|x_1-0|+w_3|x_1-0|+w_4|x_1-2|. $$
The relational table is:

\begin{center} \begin{tabular}{|l|l|l|l|l|}
\multicolumn{5}{r} {Table 5. }\\ \hline
coordinates $(a_1^{i'})$ & 0 & 0 & 0 & 2 \\ \hline
 weights $(w_{i'})$ & 4 & 1 & 2 & 3\\ \hline
$k$ & 1 & 2 & 3 &4 \\ \hline
$S_1[k]=\sum_{i=1}^{k} w_{i'}$ & 4 & 5 & 7 & 10
\\\hline
\end{tabular}
\end{center}

\smallskip
\noindent Since the equality of the form \eqref{druga}
are satisfied for $k'=2$, the solution $x_1$ is from the interval
$[a_1^{2'},a_1^{3'}]=[0,0]$, i.e. $x_1=0$, which implies $X_3(0,2)$.

\smallskip
The list $\mathfrak{X}$ is expanded: $\mathfrak{X}=\{X_1,X_2,X_3\}$.

\smallskip
Let us observe that Step $I$ transforms Table 5 into the next table

\begin{center} \begin{tabular}{|l|l|l|}
\multicolumn{3}{r} {Table 6. }\\ \hline
coordinates $(a_1^{i'})$ &  0 & 2 \\ \hline
 weights $(w_{i'})$ & 5 & 3\\ \hline
$k$ & 1 & 2  \\ \hline
$S_1[k]=\sum_{i=1}^{k} w_{i'}$ & 7 & 10 \\\hline
\end{tabular}
\end{center}
Now, condition \eqref{prva} is satisfied for $k=1$, so that $X_3(0,2)$ is possible optimal point.

\bigskip {\bf 4.} Let be $x_2\neq 1,2,4$. In this case we take $x_1=0$ and then seek for the minimum of the function
$$f_2(x_2)=\sum_{i=1}^4 w_i|x_2-a_2^i|$$

Let sort the coordinates $a_2^i\to a_2^{i'}$ and perform analogous rearrangement $w_i\to w_{i'}$. The relational table is

\begin{center} \begin{tabular}{|l|l|l|l|l|}
\multicolumn{5}{r} {Table 7. }\\ \hline
coordinates $(a_2^{i'})$ & 1 & 2 & 4 & 4 \\ \hline
 weights $(w_{i'})$ & 1 & 3 & 4 & 2\\ \hline
$k$ & 1 & 2 & 3 &4 \\ \hline
$S_2[k]=\sum_{i=1}^kw_{i'}$ & 1 & 4 & 8 & 10 \\\hline
\end{tabular}
\end{center}
\smallskip

\noindent Inequalities of the form \eqref{prva1}
hold for $k^*=3$. We stop algorithm. This case has no solution, since the assumption $x_2\neq 4$ is made.

\smallskip
Let us mention that Step $I$ gives the next Table 8.
\begin{center} \begin{tabular}{|l|l|l|l|}
\multicolumn{4}{r} {Table 8. }\\ \hline
coordinates $(a_2^{i'})$ & 1 & 2 & 4  \\ \hline
 weights $(w_{i'})$ & 1 & 3 & 6 \\ \hline
$k$ & 1 & 2 & 3  \\ \hline
$S_2[k]=\sum_{i=1}^k w_{i'}$ & 1 & 4 & 10 \\ \hline
\end{tabular}
\end{center}
Therefore, the conclusion is the same as from Table 7.

\smallskip
At the end, in order to solve Step 2 of Algorithm \ref{alg1}, we must compute and compare the values of the function $f$ at each point
$X_i$, $i=1,2,3.$ We get
$$f(X_1)=50,\mbox{  }f(X_2)=55,\mbox{  }f(X_3)=62.$$
Therefore, the solution of the Weber problem will be the point in which the function $f$ has a minimal value, i.e. $X^*=X_1=A_1=(4,4)$.
\end{exm}
\end{footnotesize}

\section {Conclusion and future work}

Our paper is the first attempt to solve the discrete and the single-facility min-sum continuous location problem
with the lift metric as the measure of distances.

\smallskip
A couple of variants and extensions of continuous location problems have been investigated in literature.
Let us mention main between them.
More complex problems include the placement of multiple facilities.
Problems with barriers are the subject in \cite{Dearing,Hamacher,Kafer,Klamroth}.
The location of undesirable (obnoxious) facilities requires to maximize
minimum distances (see, e.g., \cite{Brimberg,Erkut,[4],Melachrinoudis}.
Location models with both desirable and undesirable facilities have been analyzed in \cite{Chen1}.
It seems interesting to investigate these extensions in the sense of the lift metric or
in the more general nonconvex case, where the shortest length of arc is used as distance instead of a particular metrics.

\begin{footnotesize}

\end{footnotesize}

\end{document}